\documentclass[11pt, oneside]{amsart}
 \usepackage[text={5.58in,8.5in},centering,letterpaper,dvips]{geometry}
\usepackage{graphicx}
\usepackage{amsfonts}
\usepackage{epsf}
\usepackage{amssymb}
\usepackage{amsmath}
\usepackage{amscd}
\usepackage{tikz}
\usepackage{pdfpages}
\usepackage{fancyhdr}
\usepackage{setspace}
\usepackage{hyperref}
\usepackage[all]{xy}
\usetikzlibrary{matrix}
\usepackage{verbatim}
\usepackage{enumerate}
\usepackage{amsbsy}

\theoremstyle{theorem}
\newtheorem{theorem}{Theorem}[section]

\newtheorem{lemma}[theorem]{Lemma}

\theoremstyle{definition}

\newcommand{\Z}{\mathbb{Z}}

\newcommand{\N}{\mathbb{N}}

\newcommand{\R}{\mathbb{R}}
\newcommand{\RP}{\mathbb{RP}}

\makeatletter
\def\@seccntformat#1{%
  \protect\textup{\protect\@secnumfont
    \ifnum\pdfstrcmp{subsection}{#1}=0 \bfseries\fi
    \csname the#1\endcsname
    \protect\@secnumpunct
  }%
}  
\makeatother

\makeatletter
\newtheorem*{rep@theorem}{\rep@title}
\newcommand{\newreptheorem}[2]{%
\newenvironment{rep#1}[1]{%
 \def\rep@title{#2 \ref{##1}}%
 \begin{rep@theorem}}%
 {\end{rep@theorem}}}
\makeatother

\newreptheorem{theorem}{Theorem}
\newreptheorem{lemma}{Lemma}
\newreptheorem{question}{Question}
\newreptheorem{corollary}{Corollary}
\newreptheorem{proposition}{Proposition}
\newreptheorem{problem}{Problem}

\begin{document}

\rhead{\thepage}
\lhead{\author}
\thispagestyle{empty}


\raggedbottom
\pagenumbering{arabic}
\setcounter{section}{0}


\title{ On 2-knots and connected sums with projective planes}
\date{\today}

\author{Vince Longo}
\address{Department of Mathematics, University of Nebraska-Lincoln, Lincoln, NE 68588}
\email{vlongo2@unl.edu}
\urladdr{https://www.math.unl.edu/~vlongo2/} 
%
%

\begin{abstract}
In this paper, we generalize a result of Satoh to show that for any odd natural $n$, the connected sum of the $n$-twist spun sphere of a knot $K$ and an unknotted projective plane in the 4-sphere is equivalent to the same unknotted projective plane. We additionally provide a fix to a small error in Satoh's proof of the case that $K$ is a 2-bridge knot.
\end{abstract}

\maketitle

\section{Introduction}\label{sec:outline}
One of the earliest examples of knotted surfaces in $S^4$ is Artin's spun knot (\cite{Artin}), later generalized to twist spun knots by Zeeman in \cite{Zeeman}. These knots have easily computable fundamental groups and canonical broken surface diagrams associated to them (see \cite{Satoh2}); as such, they both provide a good starting point for many interesting questions about knotted surfaces. Problem 4.58 on Kirby's list asks whether the connected sum of an unknotted projective plane with an odd twist spun knot is always equivalent (via diffeomorphism) to the unknotted projective plane; these (connected sums of) knots all have associated fundamental group $\Z_2$ (although in the even case, their groups are typically not cyclic). By Freedman's topological s-cobordism theorem, it follows that there is a pairwise homeomorphism between $(S^4,\RP^2)$ and $(S^4,K\#\RP^2)$, where $K$ is an odd twist spun knot. By Theorem 1 of \cite{Baykur}, it follows that these knots become smoothly isotopic after enough internal stabilizations. While it is still unknown whether the knots are diffeomorphic without stabilizing, in this paper we generalize a result of Satoh to show that they become smoothly isotopic after a single (trivial) internal stabilization for any knot $K$. This falls in line with each of the examples given in \cite{Baykur}, where one internal stabilization is all that is needed to make a pair of exotically embedded surfaces smoothly isotopic. We invite the reader to compare these results with those in \cite{Hannah}, \cite{MR1953330}, and \cite{Baykur2} where one \emph{external} stabilization is all that is needed to make certain exotic smooth structures on 4-manifolds diffeomorphic.

%
%
%


%
%
%
%

%
%
 \section{Preliminaries}
 Let $K^+$ be a knotted arc in $\R_+^3=\{(x,y,z)\in \R^3| z\ge 0\}$ with endpoints in the boundary. The spin of $K$, defined originally in \cite{Artin}, is obtained by rotating the pair $(\R^3_+, K^+)$ around the $R^2$ axis given by $z=0$ and taking the continuous trace of $K^+$. If additionally we twist $K^+$ a full $n$ times while we spin around the axis (where we imagine the ``knotted" portion of $K^+$ to be contained in a rotating 3-ball), we obtain the $n$-twist spun $K$ sphere, denoted $\tau^nK$. The $n$-twist spun $K$ sphere was originally defined by Zeeman in ~\cite{Zeeman} (cf. ~\cite{Boyle} ~\cite{Surfacesbook} ~\cite{Satoh2}). If instead we replace $K^+$ with a  knot $K$ that does not intersect the boundary of $\R^3_+$, we obtain the $n$-twist spun $K$ torus, which we denote by $\sigma^nK$ to be consistent with the notation used in  ~\cite{Satoh} (cf. ~\cite{Boyle}). If $h$ is a 1-handle whose core is contained in the axis of twisting of $\tau^nK$, then $\tau^nK+h\cong \sigma^nK$, where $\tau^nK+h$ denotes performing surgery along the 1-handle $h$. Let $P_g(e)$ be an unknotted and non-orientable surface knot in $S^4$ specified by its genus $g$ and Euler number $e$ (cf. ~\cite{Kamada}). The main theorem of this paper is the following:
 \begin{theorem}
 Let $K$ be a classical knot and $n$ be a natural number. If either $n$ is odd or $K$ is a 2-bridge knot, then
 $$\tau^nK\#P_3(\pm2)\cong \tau^{n+2}K\#P_3(\pm2).$$
 \end{theorem}
 Note that this theorem implies that for all odd $n$, $\tau^nK\#P_3(\pm2)\cong P_3(\pm2)$, since $\tau^1K$ is an unknotted sphere for any knot $K$. A proof of this theorem in the case that $K$ is a 2-bridge knot and $n$ is any natural number is presented in ~\cite{Satoh}, of which the proof of our theorem heavily draws from; however, this proof has a minor error, which we will point out and show how to fix.  
\section{Proof of the Main Theorem} 
 The following lemma is proved in ~\cite{Satoh}. 
 
 \begin{lemma}\label{lem1}
 For any classical knot $K$ and any natural number $n$, $$\sigma^nK\#P_1(\pm2)\cong\sigma^{n+2}K\#P_1(\pm2).$$
 \end{lemma}
  
We present the following lemma.

\begin{lemma}\label{lem2}
Let $h$ be the 1-handle attached to $\tau^nK\#P_1(\pm2)$ whose core is contained in the axis of twisting of $\tau^nK$. Suppose that either $n$ is odd or $K$ is a 2-bridge knot. Then $h$ is isotopic to a trivial 1-handle. 
\end{lemma}
  
  The proof in ~\cite{Satoh} states that if $h$ is merely attached to $\tau^nK$ (instead of $\tau^nK\#P_1(\pm2)$), then according to ~\cite{Boyle}, $h$ is trivial (for K a 2-bridge knot and $n$ arbitrary); however, this is false in general. While Theorem 14 of ~\cite{Boyle} did show that $h$ is trivial for $n=1,2$, Theorem 15 showed that if $K$ is the trefoil or figure-8 knot, then $h$ is nontrivial for all $n\neq1,2$.

  
  \begin{proof}
Let $K$ be a classical knot in $S^3$ and let $\langle A | R \rangle$ be a Wirtinger presentation for the knot group $\pi_1(S^3\setminus K)$. Write $A=\{a_1, \ldots, a_k\}$.  In ~\cite{Satoh2}, it is shown that $$\langle A | R \cup \{a_1^na_ia_1^{-n}=a_i| i=2, \ldots k\} \rangle$$  is a presentation for the surface knot group $\pi_1(S^4\setminus \tau^nK)$. Note that connected summing a surface knot $S$ with $P_1(\pm2)$ results in giving the relation $a^2=1$ for some meridional generator $a$ of $\pi_1(S^4\setminus S)$. Thus,  $$\langle A | R\cup \{a_1^2=1\} \cup \{a_1^na_ia_1^{-n}=a_i| i=2, \ldots k\} \rangle$$ is a presentation for $\pi_1(S^4\setminus (\tau^nK\#P_1(\pm2)))$. 

Depending on the parity of $n$, this presentation can be simplified; since $a_1^2=1$, then the relation $a_1^na_ia_1^{-n}=a_i$ is equivalent to $a_1a_ia_1^{-1}=a_i$ for $n$ odd, and equivalent to the trivial relation $a_i=a_i$ for $n$ even. For $n$ odd, we see that this presentation is equivalent to the presentation $$\langle A | R \cup \{a_1^2=1\} \cup \{a_1a_i=a_ia_1 | i=2, \ldots k\}\rangle.$$ Since all of the generators are conjugate to $a_1$, this presentation is equivalent to $$\langle a_1 | a_1^2=1\rangle\cong \Z_2.$$ It was proved in ~\cite{Baykur} that this implies every handle attached to $\tau^nK\#P_1(\pm2)$ is trivial. 
For $n$ even, the presentation is equivalent to the presentation $$\langle A | R \cup \{a_1^2=1\}\rangle\cong \pi_1(S^3\setminus K)/<\mu^2>^N$$ where $\mu$ is a meridional generator for $\pi_1(S^3\setminus K)$ and $<\mu^2>^N$ is the normal closure of the subgroup generated by $\mu^2$.

For each $n$, let $\lambda_n$ denote the image of a preferred longitude $\lambda'\in \pi_1(S^3\setminus K)$ of the knot $K$ under the inclusion $(\iota_{n})_*:(B^3\setminus K^+)\times 0\to S^4\setminus \tau^nK $, where $(B^3,K^+)$ is the 3-ball, knotted arc pair (here we are using the equivalent definition of $\tau^nK$ as given in ~\cite{Boyle}). Additionally, for each $n\in \N$, write $G_{1,n} = \pi_1(S^4\setminus \tau^nK)$ and $G_{2,n}=\pi_1(S^4\setminus \tau^nK\#P_1(\pm2)).$ Note that $\langle A | R \cup \{a_1^2=1\}\rangle$ is a group presentation for each $G_{2,n}$ (and hence they are all isomorphic); as such, we will instead write $G_2$ for each $G_{2,n}$. By Theorem 14 of ~\cite{Boyle}, $\lambda_2$ is trivial. Now, each $\lambda_n$ can be represented by the same word $w\in A^*$ (specifically, $\lambda' \in\pi_1(S^3\setminus K)$ can be represented by a word in $w\in  A^*$. Now take quotients to get that $\bar{w}=_{G_{1,n}}\lambda_n$, where $\bar{w}$ is the group element that the word $w$ represents in $G_{1,n}$). Furthermore, $w$ also represents the image $\lambda_n'$ of each $\lambda_n$ under the quotient map $q_n:G_{1,n}\to G_{2,n}=G_2$ that sends $a_1^2$ to the identity. The relations for these groups are all equivalent and hence $\lambda_n'=_{G_{2}}\lambda_m'$ for all even $m$ and $n$. Since $\lambda_2$ is trivial, then $\lambda_n'=_{G_2}1$ for all even $n$. 

Let $P^+$ be the positive peripheral subgroup of $\pi_1(S^4\setminus \tau^nK\#P_1(\pm2))$. It was proved in ~\cite{Kamada} that two one handles with oriented cores $(B,C)$ and $(B',C')$ attached to a non-orientable surface knot are equivalent if and only if $P^+(B,C)P^+=P^+(B',C')P^+$. If we drag the 1-handle $h$ with its core $c$ (arbitrarily oriented) along the knot $K$ in $(B^3,K)\times \{0\} \subset \tau^nK$ so that both of the basepoints of its core are near the south pole, we can see that the core of $h$ is equivalent to the longitude $\lambda$ in $\pi_1(S^3\setminus K)$ or its inverse. Since the image of $\lambda$ (under the above maps) in $\pi_1(S^4\setminus \tau^nK\#P_1(\pm2))$ is trivial, we see that $P^+(h,c)P^+=P^+$ and hence $h$ is trivial. 

%
  
  \end{proof}
  
  With the hard part out of the way, we now prove the main theorem. This proof is essentially the same as the one given in ~\cite{Satoh}. 
  
  \begin{proof}
 Assume either $n$ is odd or $K$ is a 2-bridge knot. By Lemma \ref{lem2}, $$\tau^nK\#P_3(\pm2)\cong(\tau^nK\#P_1(\pm2))+h\cong\sigma^nK\#P_1(\pm2).$$ Then by Lemma \ref{lem1}, $$\sigma^nK\#P_1(\pm2)\cong\sigma^{n+2}K\#P_1(\pm2)\cong\tau^{n+2}K\#P_3(\pm2).$$
  \end{proof}

  \section{Acknowledgements}
  The author would like to thank his advisors Alex Zupan and Mark Brittenham for all their help, suggestions, and support throughout the process of writing this paper.

\bibliographystyle{amsalpha}
\bibliography{vince}

 \end{document}